\documentclass[12pt]{amsart}
\usepackage{amssymb}
\usepackage{amscd}
\usepackage{amsmath}

\newcommand{\nc}{\newcommand}
\nc{\cal}{\mathcal} 
\nc{\la}{\langle} \nc{\ra}{\rangle}

 \nc{\CA}{\cal A} \nc{\CBB}{\cal B}

 \nc{\CC}{\cal C}

 \nc{\CDD}{\cal D}
\nc{\CE}{\cal E}
\nc{\CF}{\cal F}
\nc{\CG}{\cal G}
\nc{\CH}{\cal H}
\nc{\CI}{\cal I}
\nc{\CJ}{\cal J}
\nc{\CK}{\cal K}
\nc{\CL}{\cal L}
\nc{\CM}{\cal M}
\nc{\CN}{\cal N}
\nc{\CO}{\cal O}
\nc{\CP}{\cal P}
\nc{\CQ}{\cal Q}
\nc{\CR}{\cal R}
\nc{\CS}{\cal S}
\nc{\CT}{\cal T}
\nc{\CU}{\cal U}
\nc{\CV}{\cal V}
\nc{\CW}{\cal W}
\nc{\CZ}{\cal Z}

\nc{\fa}{\mathfrak a}
\nc{\fg}{\mathfrak g}
\nc{\fk}{\mathfrak k}
\nc{\fh}{\mathfrak h}
\nc{\fm}{\mathfrak m}
\nc{\fn}{\mathfrak n}
\nc{\fA}{\mathfrak A}
\nc{\fC}{\mathfrak C}
\nc{\fI}{\mathfrak I}
\nc{\fL}{\mathfrak L}
\nc{\fS}{\mathfrak S}

\nc{\nen}{\newenvironment}
\nc{\ol}{\overline}
\nc{\ul}{\underline}
\nc{\lra}{\longrightarrow}
\nc{\lla}{\longleftarrow}
\nc{\Lra}{\Longrightarrow}
\nc{\Lla}{\Longleftarrow}
\nc{\Llra}{\Longleftrightarrow}
\nc{\hra}{\hookrightarrow}
\nc{\iso}{\overset{\sim}{\lra}}
\nc{\Hom}{\mathrm{Hom}}

\nc{\notebox}[1]{\noindent\fbox{\parbox{12.5cm}{\sf #1}}\\[8pt]}
\nc{\Thm}[1]{Theorem~\ref{#1}} \nc{\Prop}[1]{Proposition~\ref{#1}}
\nc{\Lem}[1]{Lemma~\ref{#1}} \nc{\Cor}[1]{Corollary~\ref{#1}}
\nc{\Conj}[1]{Conjecture~\ref{#1}} \nc{\Claim}[1]{Claim~\ref{#1}}
\nc{\Defn}[1]{ Definition~\ref{#1}} \nc{\Exa}[1]{Example~\ref{#1}}
\nc{\Rem}[1]{Remark~\ref{#1}} \nc{\Note}[1]{Note~\ref{#1}}

\nc{\marg}{\marginpar}

\nen{thm}[1]{\label{#1}{\bf Theorem.\ } \em}{}
\nen{prop}[1]{\label{#1}{\bf Proposition.\ } \em}{}
\nen{lem}[1]{\label{#1}{\bf Lemma.\ } \em}{}
\nen{klem}[1]{\label{#1}{\bf Key Lemma.\ } \em}{}
\nen{cor}[1]{\label{#1}{\bf Corollary.\ } \em}{}
\nen{conj}[1]{\label{#1}{\bf Conjecture.\ } \em}{}
\nen{thma}[1]{\label{#1}{\bf Theorem A.} \em}{}
\nen{thmb}[1]{\label{#1}{\bf Theorem B.\ } \em}{}

\nen{claim}[1]{\label{#1}{\bf Claim.\ } \em}{}

\nen{defn}[1]{\label{#1}{\bf Definition.\ } }{}
\nen{exa}[1]{\label{#1}{\bf Example.\ } }{}

\nen{rem}[1]{\label{#1}{\em Remark.\ } }{}
\nen{exer}[1]{\label{#1}{\em Exercise.\ } }{}
\nen{pf}{\begin{proof}}{\end{proof}}

\nc{\ba}{\mathbb A}
\nc{\br}{\mathbb R} \nc{\bz}{\mathbb Z} \nc{\bc}{\mathbb C}
\nc{\bn}{\mathbb N} \nc{\geg}{\mathfrak g} \nc{\gh}{{\goth h}}
\nc{\gan}{{\goth a}} \nc{\G}{\Gamma} \nc{\g}{\gamma}
\nc{\sm}{\setminus} \nc{\sub}{\subset} \nc{\lm}{\lambda}
\nc{\eps}{\varepsilon} \nc{\nty}{\infty} \nc{\al}{\alpha}
\nc{\bt}{\beta} \nc{\om}{\omega} \nc{\dl}{\delta} \nc{\Dl}{\Delta}
\nc{\Om}{\Omega} \nc{\s}{\sigma} \nc{\ro}{\rho} \nc{\te}{\theta}
\nc{\SLR}{SL_2(\br)} \nc{\GLR}{GL_2(\br)} \nc{\PGLR}{PGL_2(\br)}
\nc{\PSLR}{PSL_2(\br)}
\nc{\PSLZ}{PSL_2(\bz)}
\nc{\SLC}{SL(2,\bc)}
\nc{\PGL}{PGL_2}
\nc{\uH}{\mathbb H}
\nc{\fD}{\mathfrak D} \nc{\fE}{\mathfrak E} \nc{\haf}{\frac{1}{2}}
\nc{\qtr}{\frac{1}{4}} \nc{\8}{\infty} \nc{\7}{{-\infty}}
\nc{\inv}{^{-1}}
\nc{\shaf}{{\scriptstyle\frac{1}{2}}}
\nc{\hlm}{{\scriptstyle\frac{\lambda}{2}}}
\begin{document}


\title[automorphic functions]
{Estimates of triple products of automorphic functions II }

\author{Andre Reznikov}
\address{Bar Ilan University, Ramat-Gan, Israel}
\email{reznikov@math.biu.ac.il}

\begin{abstract}
We prove a sharp bound for the average value of  the triple product of modular functions for the Hecke subgroup $\G_0(N)$.  Our result is an extension of the main result in \cite{BR-moscj} to a {\it fixed} cuspidal representation of the adele group $\PGL(\ba)$.
\end{abstract}
\thanks{Partially supported by  the ERC grant, by the ISF  Center of Excellency grant 1691/10, and by the Minerva  Center at ENI}
\maketitle

\section{Introduction}
\label{intro}
\subsection{Maass forms} We recall the setup of \cite{BR-moscj} which should be read in conjunction with this appendix.
Let $Y$ be a compact Riemann surface  with a Riemannian metric of
constant curvature $-1$ and the associated volume element $dv$.
The corresponding  Laplace-Beltrami operator  is non-negative and
has purely discrete spectrum on the space $L^2(Y,dv)$ of functions
on $Y$. We will denote by $0=\mu_0< \mu_1 \leq \mu_2 \leq ...$ its
eigenvalues  and by $\phi_i$ the corresponding
eigenfunctions (normalized to have $L^2$ norm one). In the theory
of automorphic forms the functions $\phi_i$ are called
automorphic functions or Maass forms (after H. Maass,
\cite{M}).  We write $\mu_i=(1-\lm_i^2)/2$  and $\phi_i=\phi_{\lm_i}$ as  is customary in representation theory of the group
$\PGLR$.

For any three Maass forms
$\phi_i,\ \phi_j, \ \phi_k$ we define the following
triple product or triple period:
\begin{eqnarray}\label{cijk}
c_{ijk}=\int_Y\phi_i\phi_j\phi_kdv \ .
\end{eqnarray}

One would like to bound the coefficient $c_{ijk}$ as a function of
eigenvalues  $\mu_i,\ \mu_j,\ \mu_k$. In particular, we would like
to find bounds for these coefficients  when one or more of these
indices tend to infinity. The study of these triple coefficients goes back  to
pioneering works of Rankin and Selberg (see \cite{Ra}, \cite{Se}), and reappeared in celebrated works of Waldspurger \cite{W} and Jacquet \cite{J} (also see  \cite{HK}, \cite{Wa}, \cite{Ich}). Recently, an interest in analytic questions related to triple products was initiated in the groundbreaking paper of Sarnak  \cite{Sa} (see also  \cite{Go} for the first general result on exponential decay). This was motivated by the widespread use of triple products in applications  (e.g., \cite{Sa2}).

In our paper \cite{BR-moscj} we considered the following problem. We fix two
Maass forms $\phi = \phi_\tau,\ \phi'= \phi_{\tau'}$ as above and
consider coefficients defined by the triple period:
\begin{eqnarray}\label{ci}
c_i=\int_Y\phi\phi'\phi_idv
\end{eqnarray}
as  $\{\phi_i = \phi_{\lm_i}\}$ run over the orthonormal basis of Maass forms.
We note first that one  has exponential decay for the coefficients $c_i$
in the parameter $|\lm_i|$ as $i$ goes to $\8$. For that reason, one renormalizes coefficients $|c_i|^2$ by an appropriate ratio of Gamma functions dictated by the Watson formula \cite{Wa} (see also Appendix in \cite{BR-moscj} where these factors were computed from another point of view). Taking into account the asymptotic behavior of these factors, we introduced normalized coefficients
\begin{eqnarray}\label{ai-ci}
|a_i|^2=|\lm_i|^2\exp\left(\frac{\pi}{2}|\lm_i|\right)\cdot |c_i|^2\ .
\end{eqnarray} Under such a normalization, we showed that
\begin{eqnarray}\label{mean-value-ai}
\sum_{ |\lm_i|\leq T}|a_i|^2\leq A\cdot T^2\
\end{eqnarray} for some explicit constant $A=A(\G,\phi,\phi')$.
According to the Weyl law, there are approximately $cT^2$ terms in the above sum, and hence the bound \eqref{mean-value-ai} consistent with the Lindel\"{o}ff  bound on average (in fact it is not difficult to show that the bound \eqref{mean-value-ai} essentially is  sharp, see \cite{Re1}).

There are various natural questions concerning the bound \eqref{mean-value-ai} which were not discussed in \cite{BR-moscj}. These are mostly related to the dependence of the constant $A$ on various parameters (i.e., $\G$, $\phi$, $\phi'$), and also to the fact that we restricted the discussion to Maass forms, leaving aside the case of holomorphic forms. Another restriction of  the treatment we presented was  the fact that we used in an essential way the compactness of $Y$ (which was treated by us differently in \cite{BR1}).  All these issues turn out to be important in applications. In this appendix, we answer some of these questions for Hecke congruence subgroups.
\subsubsection{Hecke subgroups} We consider Hecke subgroups $\G_0(N)$ of the modular group. We normalize the scalar product on the quotient Riemann surface $Y_{N}=\G_0(N)\setminus \uH$ by $\langle f,g\rangle_{Y_N}=
\frac{1}{{\rm vol}_\uH(Y_N)}\int_{Y_{N}}f(z)\bar g(z)d\mu_\uH$  where $d\mu_\uH$ is the standard volume element on the upper half plane $\uH$ (i.e., we normalize the volume element $dv_{Y_{N}}$ on $Y_N$ to have the total volume $1$). Let $\phi$ be a (primitive) Hecke-Maass form for the group $\G_0(N_0)$ for some fixed level $N_0\geq 1$. We assume that $\phi$ is normalized by the $L^2$-norm $||\phi||_{L^2(Y_{N_0})}=1$.  For an integer $\ell> 1$, denote by $^\ell\phi(z)=\phi(\ell z)$  the corresponding {\it old form} for the Hecke subgroup $\G_0(\ell N_0)$. The corresponding function $^\ell\phi$  also turns out  to be  $L^2$-normalized on $Y_{\ell N_0}$ with respect to the normalization of measures we choose. This follows easily from the Rankin-Selberg method. For two such Maass forms $\phi$ and $\phi'$, we define  triple products by
\begin{eqnarray}\label{ci-ell}
c_i(\ell)=\int_{Y_{\ell N_0}}{^\ell\phi}\ {^\ell\phi'}\ \phi_i\ dv_{Y_{\ell N_0}}
\end{eqnarray}
as  $\{\phi_i = \phi_{\lm_i}\}$ run over the orthonormal basis of Maass forms on $Y_{\ell N_0}$. We have the corresponding  normalized coefficients
\begin{eqnarray}\label{ai-ell}
|a_i(\ell)|^2=|\lm_i|^2\exp\left(\frac{\pi}{2}|\lm_i|\right)|c_i(\ell)|^2\ .
\end{eqnarray}

\begin{thma}{THM1}
There exists an effectively computable
constant $A$  such that the following bound holds for all $T \geq 1$,
\begin{eqnarray}\label{mean-value}
\sum_{ |\lm_i|\leq T}|a_i(\ell)|^2\leq A\cdot T^2\ ,
\end{eqnarray} where the summation is over an orthonormal basis of Maass forms for the group $\G_0(\ell N_0)$.
The constant $A$ depends on $N_0$, $\phi$ and $\phi'$, but not on $\ell$.
\end{thma}

This could be viewed as the Lindel\"{o}ff on the average bound in two parameters $\ell$ and $T$ (but only for old forms of course). There are $[\G_0(\ell N_0):\G(1)]\cdot  T^2$ terms in the sum according to the Weyl-Selberg law, and the resulting bound is consistent with this under our normalization of measures   ${\rm vol}(Y_{\ell N_0})=1$. In fact, there is also an analogous contribution from the Eisenstein series leading to the triple coefficients $a_{s,\kappa}(\ell)$ arising from triple products $\langle {^\ell\phi}\ {^\ell\phi'}, E_{\kappa}(s)\rangle_{Y_{\ell N_0}}$ where $E_{\kappa}(s)$ is the Eisenstein series associated to a cusp $\kappa$ of $\G_0(\ell N_0)$. We have then
\begin{eqnarray}\label{mean-value-Eisen}
\sum_{ |\lm_i|\leq T}|a_i(\ell)|^2+\sum_{\{\kappa\}}\int_{-T}^{T}|a_{it,\kappa}(\ell)|^2dt \leq A\cdot T^2\ .
\end{eqnarray}

\subsection {The method}
The proof we presented in \cite{BR-moscj} was based on the uniqueness of the triple
product in representation theory of the group $\PGLR$.
We review quickly the general ideas behind our  proof. It is
based on ideas from representation theory (see the seminal book \cite{GGPS}, and also \cite{Bu}, \cite{L}). Namely, we use the fact
that every automorphic form $\phi$ generates an automorphic
representation of the group $G = \PGLR$; this means that starting
from $\phi$ we produce a smooth irreducible representation of the
group $G$ in a space $V$ and its realization $\nu : V \to
C^{\8}(X)$ in the space of smooth functions on the automorphic
space $X = \G \backslash G$ endowed with the invariant measure of the total mass one. We denote by $\phi_v(x)=\nu(v)(x)$ the corresponding to $v\in V$ automorphic function. The Maass form corresponds to a unit $K$-invariant vector $e_0\in V$, where $K=PSO (2)$ is a maximal compact connected subgroup of $G$. In our case, we have the family of spaces $X_\ell=\G_0(\ell N_0)\setminus G$ (we denote by $X_0=\G_0( N_0)\setminus G$), and the the corresponding family of  isometries $\nu_\ell: V\to C^{\8}(X_\ell)$ of the {\it same} abstract representation of $G$. These are generated by old Maass forms $\phi(\ell z)$ in the same adelic representation corresponding to the Hecke-Maass form $\phi$.

 The triple product  $c_i=\int_Y\phi\phi'\phi_idv$ extends to a
$G$-equivariant trilinear form on the corresponding automorphic
representations  $l^{aut}:V\otimes V'\otimes V_i\to\bc$, where $V
= V_\tau, V' = V_{\tau'}, V_i = V_{\lambda_i}$ .

Then we use a general result from representation theory that such a
$G$-equivariant trilinear form is unique up to a scalar, i.e., that the space $\Hom_G(V\otimes V'\otimes V_i,\bc)$ is at most one-dimensional (see \cite{O},\cite{Mo}, \cite{Lo} and \cite{Pr} for $p$-adic $GL(2)$).
This
implies that the automorphic form $l^{aut}$ is proportional to an
explicit ``model" form $l^{mod}$ which we describe using explicit
realizations of representations of the group $G$; it is important
that this last form carries no arithmetic information.

Thus we can write $l^{aut} = a_i \cdot l_i^{mod}$ for some constant
$a_i$ and hence $c_i= l_i^{aut}(e_\tau\otimes e_{\tau'} \otimes
e_{\lambda_i}) = a_i \cdot l_i^{mod}(e_\tau \otimes e_{\tau'}
\otimes e_{\lambda_i})$, where $e_\tau,\ e_{\tau'},\
e_{\lambda_i}$ are $K$-invariant unit vectors in the automorphic
representations $V, V',\ V_i$ corresponding to the automorphic
forms $\phi$, $\phi'$ and $\phi_i$.

 It turns out that the proportionality coefficient
 $ a_i$ in the last formula carries
 important ``automorphic" information while the second factor
 carries no arithmetic information and can be computed
 in terms of Euler $\G$-functions using explicit realizations of
 representations $V_\tau$, $V_{\tau'}$ and $V_{\lambda_i}$ (see Appendix to \cite{BR-moscj}). This second factor is
 responsible for
 the exponential decay, while the first factor $a_i$ has a
 polynomial behavior in the parameter $\lm_i$.

   In order to bound the quantities $a_i$, we use the fact that they
   appear as coefficients in the spectral
   decomposition of the diagonal Hermitian form $H_{\Dl}$ given by $$H_\Dl(v\otimes w)=\int_X|\phi_v(x)\phi_w(x)|^2dx$$
   on the space $E = V_\tau \otimes V_{\tau'}$.
   This gives an inequality $\sum |a_i|^2 H_i \leq H_\Dl$
   where $H_i$ is a Hermitian
form on $E$ induced by the model trilinear form $l_i^{mod}: V
\otimes V' \otimes V_i \to \bc$ as above.

Using the geometric properties of the diagonal form and simple explicit estimates of forms $H_i$,
we establish the mean-value bound for the coefficients $|a_i|^2$.
Here is where one obtains the dependence of the constant $A$ in \eqref{mean-value-ai} on parameters involved. In the method of \cite{BR-moscj}, we used $L^2$ theory by averaging the form $H_\Dl$ and comparing the resulting form with the $L^2$-form. The coefficient $A$ that one obtains in such an argument depends in particular on the injectivity radius of $X$. While in certain cases it gives an optimal result, it obviously has two drawbacks. One is related to the possible non-compactness of $X$ since in the cusp the injectivity radius tends to zero.  Another problem arises when one considers a sequence of subgroups with co-volume going to infinity along the sequence. In that case, the bound which the method of \cite{BR-moscj} provides for the constant (e.g., $A(\G_0(p))\leq {\rm vol}(Y_p)\approx p$, see \cite{K})  is too weak for many applications. Both of these problems arise in the classical setup of Hecke subgroups $\G_0(N)$. Here we obtain an optimal bound for {\it old forms}. We do not know how to obtain similar results for new forms. We will discuss improvements over a trivial bound for new forms elsewhere. Theorem \ref{THM1} could be viewed as the exact analog of the result in \cite{BR-moscj} for a {\it fixed} adelic representation.

\subsection{Proof of Theorem A} We have a family of various objects  $(X_\ell, \nu_\ell, H_\Dl^{\ell}, a_i(\ell))$ parameterized by the level $\ell$. However the {\it model } Hermitian form $H_i$ is the same since the abstract representation $V$ of $\PGLR$ does not change.
The proof of the bound \eqref{mean-value-ai} given in \cite{BR-moscj} was based on the spectral decomposition $\sum |a_i(\ell)|^2 H_i \leq H^\ell_\Dl$ of the diagonal form and on the construction of the test vector $u=u_T\in E$ such that $H^\ell_\Dl(u)\leq aT^2$ and $H_i(u)\geq 1$ for $|\lm_i|\leq T$. We now construct the test vector {\it independently} of $\ell$. The dependence on $\ell$ is hidden in the automorphic realization $\nu_\ell(u)$ as a function on $X_\ell$.
\subsubsection{Construction of vector $u$}\label{consttest} We slightly change the construction of the test vector $u_T$ given in Section 5.3.2 of \cite{BR-moscj}.
Let us identify the space $E = V \otimes V'$ with a  subspace of smooth functions $C^\8(\br\times\br)$. Choose a smooth non-negative function $\al\in C^\8(\br)$ with the support ${\rm supp}(\al)\subset [-0.1,0.1]$ and $\int_\br\al(x)dx=1$. Let $||\al||^2_{L^2(\br)}=c^2$ for some $c>0$. Consider the diagonal element $a_T={\rm diag}(T^{-\haf},T^\haf)\in G$. We define  vectors
\begin{eqnarray}\label{test-vect}
v_T=T^{\haf+\tau}\cdot\pi_\tau(a_T)\al\ \ {\rm and}\ \  v'_T=T^{\haf+\tau'}\cdot\pi_{\tau'}(\left(
                                          \begin{smallmatrix}
1 & -1 \\
                                   &\ \ 1 \\
\end{smallmatrix}
                                        \right)) v_T\ .
                                        \end{eqnarray}
We set our test vector to $u_T(x,y)= v_T(x)\otimes v'_T(y)$.  

Recall that the action is given by $\pi_\tau({\rm diag}(a\inv,a) v(x)=|a|^{1-\tau}v(a^2x)$, and $\pi_{\tau}(\left(
                                          \begin{smallmatrix}
1 & -1 \\
                                   &\ \ 1 \\
\end{smallmatrix}
                                        \right)) v(x)=v(x+1)$.  By an easy calculation we have $||u_T||^2_E=c^2T^2$. We note that geometrically the vector $u_T$ is a small non-negative bump function around the point $(0,1)\in \br^2$, with the support in the box of the size $T\inv$, and satisfies $\int_{\br^2}u_T(x,y)dxdy=1$.
Computation identical to the one performed in Section 5.3.4 of \cite{BR-moscj} gives then $H_i(u_T)\geq \beta$ for $|\lm_i|\leq T$ and some explicit $\beta>0$ independent of $\lm_i$ (in fact for $|\tau|,\ |\tau'|\leq T$, it does not depend on these either). We remark that the only difference with the construction of the test vector given in \cite{BR-moscj} is that here we constructed $u_T$ with the help of the action of $G$ on $V$ (while in \cite{BR-moscj} we constructed essentially the same vector explicitly in the model). This will play a crucial role in our estimate of the corresponding automorphic function.

We now need to estimate $H^\ell_\Dl(u_T)$.  We claim that $H^\ell_{\Dl}(u_T)\leq B T^2$ for some explicit constant $B$ independent of $\ell$. Since we have
$$H^\ell_{\Dl}(u_T)=\int_{X_\ell}|{^\ell\phi_{v_T}}{^\ell\phi_{v'_T}}|^2dx\leq \haf||{^\ell\phi_{v_T}^2}||^2_{L^2(X_\ell)} +\haf||{^\ell\phi_{v'_T}^2}||^2 _{L^2(X_\ell)}\ ,$$ it is enough to show that $||{^\ell\phi_{v_T}^2}||^2_{L^2(X_\ell)}\leq \beta' T^2$.  This would finish the proof of Theorem~\ref{THM1} following the argument  in Section 4.7 \cite{BR-moscj}.  Since $\langle {^\ell\phi_{v_T}^2}, 1\rangle_{X_\ell}=||{^\ell\phi_{v_T}}||^2=c T$, it is easy to see that such a bound is sharp. We claim that
\begin{eqnarray}\label{supnorm}
\sup_{x\in X_\ell}|{^\ell\phi_{v_T}}(x)|\leq \beta'' T^\haf\ ,
\end{eqnarray} for some $\beta''$ independent of $\ell$. Here  ${^\ell\phi_{v_T}}(x)=\nu_\ell(v_T)(x)$. Note that the bound provided by the Sobolev theorem \cite{BR2} is too weak.

\subsubsection{Supremum norm}\label{supnorm-sect} Recall that we started with an $L^2$-normalized Hecke-Maass form $\phi$ on $\G_0(N_0)$, and the corresponding isometry of $\nu:V\to C^\8(X_0)$ of the principal series representation $V\simeq V_\tau$. We then constructed another isometry $\nu^\ell:V\to C^\8(X_{\ell})$ by using the map
\begin{eqnarray}\label{nu-ell-nu}
^\ell\phi_v(x)=\nu^\ell(v)(x)=\\ \nu(v)\left(\left(
                                          \begin{smallmatrix}
\ell^\haf &  \\
                                   &\ \ \ell^{-\haf} \\
\end{smallmatrix}
                                        \right)x\right)=\phi_v\left(\left(
                                          \begin{smallmatrix}
\ell^\haf &  \\
                                   &\ \ \ell^{-\haf} \\
\end{smallmatrix}
                                        \right)x\right)\ ,\nonumber
\end{eqnarray}        for any $v\in V$.    This relation might be viewed as the relation between  functions on $G$ invariant on the left for an appropriate $\G$ (e.g., for  $\G_0(N_0)$ and for $\G_0(\ell N_0)$). In particular, we see that the supremum of the function $^\ell\phi_v$ on $X_\ell$ and that of the function $\phi_v$ on $X_0$  are equal for the {\it same} vector $v\in V$. Hence it is enough to show that $\sup_{x\in X_0}|{\phi_{v_T}}(x)|\leq \beta'' T^\haf$. In fact, this is obvious since  $v_T=T^{\haf+\tau}\pi(a_T)\al$ is given by the (scaled) action of $G$ on a {\it fixed} vector.  We have
\begin{eqnarray*}\label{sup-action}
\sup_{X_0}|{\phi_{v_T}(x)}|=T^{\haf}\sup_{X_0}|{\phi_{\pi(a_T)\al}(x)}|=\\
T^{\haf}\sup_{X_0}|{\phi_{\al}(x\cdot a_T)}|=
T^{\haf}\sup_{X_0}|{\phi_{\al}(x)}|= \beta'' T^\haf\ ,
\end{eqnarray*}
since $\al$ is a {\it fixed} vector in a  {\it fixed} automorphic cuspidal  representation $(\nu, V)$, and the action does not change the supremum norm. \qed

\subsubsection*{Remark}   It is easy to see that the condition that forms $^\ell\phi$ and $^\ell\phi'$ be  of the same level is not essential for the proof, as well as that these be Hecke forms. In particular, under our the normalization of measures on $X_\ell$, we see that for a vector $v\in V$, $L^2$-norms $||{^\ell\phi_v}||_{X_\ell}=||{^\ell\phi_v}||_{X_{\ell'}}$ are equal  if $\ell|\ell'$ (here we view the function ${^\ell\phi_v}$ as both $\G_0(\ell N_0)$-invariant function and as
$\G_0(\ell' N_0)$-invariant function). Hence for two Maass forms  $\phi$ and $\phi'$ on $\G_0(N_0)$, we obtain the bound:
\begin{eqnarray}\label{mean-value1}
\sum_{ |\lm_i|\leq T}|\langle {^{\ell_1}}\phi\cdot{^{\ell_2}}\phi',\phi_i\rangle_{Y_{\ell_1\ell_2 N_0}}|^2\cdot|\lm_i|^2e^{\frac{\pi}{2}|\lm_i|}\leq A\cdot T^2\ ,
\end{eqnarray} where the summation is over an orthonormal basis of Maass forms for the subgroup $\G_0(\ell_1\ell_2 N_0)$.

\subsection{Holomorphic forms} The approach given above is applicable to holomorphic forms as well. In principle, there are no serious changes needed as compared to the Maass forms case. The main difficulty is that   we have to fill up the gap left in  \cite{BR-moscj} concerning the model trilinear functional for discrete series representations of $\PGLR$.

Let $\phi^k$, $\phi'^k$ be (primitive) holomorphic forms of weight $k$ for the subgroup $\G_0(N_0)$. We assume these are $L^2$-normalized. For $\ell>1$, we consider (old) forms ${^\ell\phi^k}(z)=\phi^k(\ell z)$ and ${^\ell\phi'^k=\phi'^k(\ell z)}$ on $\G_0(\ell N_0)$. Under our normalization of measures for $Y_{\ell N_0}$, we have  $||{^\ell\phi^k}||_{Y_{\ell N_0}}= ||{^\ell\phi'^k}||_{Y_{\ell N_0}}=\ell^{-\frac{k}{2}}$. This follows from the Rankin-Selberg method. Hence it would have been  more natural to consider normalized forms $\phi^k|_{[a_\ell]_k}=\ell^{\frac{k}{2}}\cdot{^\ell\phi^k}$.

For a (norm one) Maass form $\phi_i$ on $\G_0(\ell N_0)$, we define the corresponding triple coefficient by
\begin{eqnarray}\label{ci-ell-hol}
c^k_i(\ell)=\int_{Y_{\ell N_0}}{^\ell\phi^k}\ \overline{^\ell\phi'^k}\ \phi_i\ y^k\ dv_{Y_{\ell N_0}}\ .
\end{eqnarray} As with Maass forms, we renormalize these coefficients in accordance with the Watson formula by introducing normalized triple product coefficients
\begin{eqnarray}\label{ai-ell-hol}
|a^k_i(\ell)|^2=|\lm_i|^{2-2k}\exp\left(\frac{\pi}{2}|\lm_i|\right)|c^k_i(\ell)|^2\ .
\end{eqnarray}

\begin{thmb}{THM2}
There exists an effectively computable
constant $A$  such that the following bound holds for all $T \geq 1$,
\begin{eqnarray}\label{mean-value-hol}
\sum_{\haf T\leq |\lm_i|\leq T}|a^k_i(\ell)|^2\leq A\cdot {\ell^{-2k}}T^{2}\ ,
\end{eqnarray} where the summation is over an orthonormal basis of Maass forms for the group $\G_0(\ell N_0)$.
The constant $A$ depends on $N_0$, $\phi$ and $\phi'$, but not on $\ell$.
\end{thmb}

\subsubsection*{Remark} The proof we give applies to a slightly more general setup of forms  of different co-prime level $\ell_1$ and $\ell_2$.   Namely, we have
\begin{eqnarray}\label{mean-value-hol-dif-level}
\sum_{ \haf T\leq |\lm_i|\leq T}|\langle {^{\ell_1}\phi^k}\ \overline{{^{\ell_2}}\phi'^k},\phi_i\rangle_{Y_{\ell_1\ell_2 N_0}}|^2\cdot(\ell_1\ell_2)^{k}|\lm_i|^{2-2k}e^{\haf\pi|\lm_i|}\leq A\cdot T^2\ ,
\end{eqnarray} for two forms $\phi^k$ and $\phi'^k$ on $\G_0(N_0)$. Breaking the interval $[1,T]$ into dyadic parts, we obtain for the full range,
\begin{eqnarray}\label{mean-value-hol-dif-level2}\\
\sum_{ |\lm_i|\leq T}|\langle {^{\ell_1}\phi^k}\ \overline{{^{\ell_2}}\phi'^k},\phi_i\rangle_{Y_{\ell_1\ell_2 N_0}}|^2\cdot(\ell_1\ell_2)^{k}|\lm_i|^{2-2k}e^{\haf\pi|\lm_i|}\leq A\cdot T^2\ln(T)\ .\nonumber
\end{eqnarray}  This is slightly weaker than \eqref{mean-value} for Maass forms.

\subsection{Proof of Theorem B} As we seen in the case of Maass forms, the proof is based on the explicit form of the trilinear functional, its value on special vectors leading to the normalization \eqref{ai-ell-hol}, and the construction of test vectors for which we can estimate supremum norm effectively. We explain below changes and additions needed in order to carry out this scheme for discrete series.

\subsubsection{Discrete series} Let $k\geq 2$ be an even integer, and $(D_k,\pi_{D_k})$ be the corresponding discrete series representation of $\PGLR$. In particular, for $m\in 2\bz$, the space of  $K$-types of weight $m$ is non-zero (and in this case is one-dimensional) if and  only if $|m|\geq k$. This defines $\pi_k$ uniquely. Under the restriction to $\PSLR$, the representation $\pi_k$ splits into two representations $(D^\pm_k,\pi^\pm_{D_k})$ of ``holomorphic" and ``anti-holomorphic" discrete series  of $\PSLR$, and  the element $\dl={\rm diag}(1,-1)$ interchanges them.

We consider two realizations of discrete series as subrepresentations and as quotients of induced representations. Consider the space  $\CH_{k-2}$  of smooth even homogeneous functions on $\br^2\setminus 0$ of  homogeneous degree $k-2$ (i.e., $f(tx)=t^{k-2}f(x)$ for any $t\in \br^\times$ and $0\not=x\in\br^2$). We have the natural action of $\GLR$ given by ${\pi}_{k-2}(g)f(x)=f(g\inv x)\cdot \det(g)^{(k-2)/2}$, which is trivial on the center and hence defines a representation $(\CH_{k-2}, \pi_{k-2})$ of $\PGLR$. There exists a unique non-trivial invariant subspace $W_{k-2}\subset \CH_{k-2}$. The space $W_{k-2}$ is finite-dimensional, $\dim W_{k-2}=k-1$, and is generated by monomials $x_1^{m}x_2^{n}$, $m+n=k-2$. The quotient space $\CH_{k-2}/ W_{k-2}$ is isomorphic to the space of smooth vectors of the discrete series representation $\pi_k$.

We also consider the dual situation. Let  $\CH_{-k}$ be the space of smooth even homogeneous functions on $\br^2\setminus 0$ of  homogeneous degree $-k$. There is a natural $\PGLR$-invariant pairing $\langle\ ,\ \rangle: \CH_{k-2}\otimes \CH_{-k}\to\bc$ given by the integration over $S^1\subset \br^2\setminus 0$. Hence $\CH_{-k}$ is the smooth dual of $\CH_{k-2}$, and vice versa. There exists a unique  non-trivial invariant subspace $D^*_k\subset \CH_{-k}$. The quotient $\CH_{-k}/D^*_k$ is isomorphic to the finite-dimensional representation $W_{k-2}$.

Of course  $D^*_k$ is isomorphic to $D_k$, but we will distinguish between two realizations of the same abstract representation  as a subrepresentation $D^*_k\subset \CH_{-k}$ and as a quotient $\CH_{k-2}\to D_k$. We denote corresponding maps by $i_k: D^*_k\subset \CH_{-k}$ and $q_k: \CH_{k-2}\to D_k$.

\subsubsection{Trilinear invariant functionals}
Let $(V_{\lm,\eps},\pi_{\lm,\eps})$ be a unitary representation of the principal series of $\PGLR$. These are parameterized by $\lm\in i\br$ and by $\eps=0, 1$ describing the action of the element $\dl$ (see \cite{Bu}). The space $\Hom_G(D_k\otimes D^*_k, V_{\lm,\eps})$ is one-dimensional. We will work with the space of invariant trilinear functionals $\Hom_G(D_k\otimes D^*_k\otimes V_{-\lm,\eps},\bc)$ instead. We construct  below a non-zero functional  $l^{ind}_{k,\lm,\eps}\in \Hom_G(\CH_{k-2}\otimes\CH_{-k}\otimes V_{-\lm,\eps},\bc)$ for induced representations (in fact, this space is also one-dimensional) by means of (analytic continuation of) the explicit kernel.  We use it to define a non-zero functional   $l^{mod}_{k,\lm,\eps}\in\Hom_G(D_k\otimes D^*_k\otimes V_{-\lm,\eps},\bc)$. What is more important, we will use $ l^{ind}_{k,\tau,\eps}$ in order to carry out our computations in a  way similar to the principal series.

Let $ l^{ind}_{k,\lm,\eps}\in \Hom_G(\CH_{k-2}\otimes\CH_{-k}\otimes V_{-\lm,\eps},\bc)$ be a non-zero  invariant functional.  Such a functional induces the corresponding functional on $\CH_{k-2}\otimes D^*_k\otimes V_{-\lm,\eps}$ since $D^*_k\subset \CH_{-k}$. Moreover, any such functional vanishes on the subspace $W_{k-2}\otimes D^*_k\otimes V_{-\lm,\eps}$ since there are no non-zero maps between $W_{k-2}\otimes D^*_k$ and $V_{\lm,\eps}$. Hence we obtain a functional  $ l^{mod}_{k,\lm,\eps}\in \Hom_G(D_k\otimes D^*_k\otimes V_{-\lm,\eps},\bc)$ on the corresponding quotient space. We denote by $T^{mod}_{k,\lm,\eps}: D_k\otimes D^*_k\to V_{\lm,\eps}$ the associated  map, and  by $H^{mod}_{k,\lm,\eps}(u)=||T^{mod}_{k,\lm,\eps}(u)||_{V_{\lm,\eps}}^2$, $u\in D_k\otimes D^*_k$ the corresponding Hermitian form.

\subsubsection{Model  functionals} We follow the construction from \cite{BR-moscj}. Denote \begin{eqnarray}\label{K-kernel2}\ \ 
K_{k,\lm}(x,y,z)=|x-y|^{\frac{-1-\lm}{2}}|x-z|^{\frac{-1+\lm}{2}-k+1}|y-z|^{\frac{-1+\lm}{2}+k-1}\ .
\end{eqnarray}
In order to construct $ l^{mod}_{k,\lm,\eps}\in \Hom_G(\CH_{k-2}\otimes\CH_{-k}\otimes V_{-\lm,\eps},\bc)$, we consider the following function  in three variables $x,\ y,\ z\in \br$
\begin{eqnarray}\label{K-circle-discr}
K_{k-2,-k,\lm,\eps}(x,y,z)=(sgn(x,y,z))^{\eps}\cdot K_{k,\lm}(x,y,z)\ ,
\end{eqnarray}  where $sgn(x_1,x_2,z_3)=\prod_{i\not= j}sgn(x_i-x_j)$ (this is an $\SLR$-invariant function on $\br^3$ distinguishing two open orbits). An analogous  expression could be written in the circle  model on the space $C^\8(S^1)$.  Viewed as a kernel, $K_{k-2,-k,\lm,\eps}$ defines an invariant non-zero functional $l^{ind}_{k,\lm,\eps}$ on the (smooth part of) the representation $\CH_{k-2}\otimes\CH_{-k}\otimes V_{-\tau,\eps}\subset C^\8(\br^3)$.
 Such a kernel should be understood in the regularized sense (e.g.,  analytically  continued following  \cite{G1}). We are interested in $\lm\in i\br$, $|\lm|\to\8$, and hence all  exponents in \eqref{K-circle-discr} are non integer. This implies that the regularized kernel does not have a pole at relevant points.

 We denoted by $l^{mod}_{k,\lm,\eps}\in \Hom_G(D_k\otimes D^*_k\otimes V_{-\lm,\eps},\bc)$ the corresponding model functional. The difference with principal series clearly lies in the fact that we only can compute the auxiliary functional $l^{ind}_{k,\lm,\eps}$. However, for $k$ fixed, it turns out that necessary computations are essentially identical to the ones we performed for the principal series in \cite{BR-moscj}.

\subsubsection{Value on $K$-types}  In order to obtain normalization \eqref{ai-ell-hol}  and to compare our model functional $l^{mod}_{k,\lm,\eps}$ to the automorphic triple product \eqref{ci-ell-hol}, we have to compute, or at least to bound,  the value $l^{mod}_{k,\lm,\eps}(e_k\otimes e_{-k}\otimes e_0)$ where $e_{\pm k}\in D_k$ are  highest$/$lowest $K$-types of  norm one, and $e_0\in V_{\lm,\eps}$ is  a  $K$-fixed vector of norm one. For Maass forms, this is done in the Appendix of \cite{BR-moscj} by explicitly calculating this value in terms of $\G$-functions.  In fact, the relevant calculation is valid for  $K$-fixed vectors for any three induced representations with generic values of parameters (i.e., those for which the final expression is well-defined).  Using the action of Lie algebra of $G$ (see \cite{Lo} for the corresponding calculation where it is used to prove uniqueness), one obtains recurrence relations between values of  the model functional on various weight vectors.  For a generic value of $\tau$, this allows one to reduce the computation of $l^{mod}_{\tau,\lm,\eps}(e_k\otimes e_{-k}\otimes e_0)$ to the value  of $l^{mod}_{\tau,\lm,\eps}(e_0\otimes e_{0}\otimes e_0)$. By analytic continuation,  this relation holds for our set of parameters corresponding to discrete series. From this, one deduces  the bound
\begin{eqnarray}\label{K-value}
|l^{mod}_{k,\lm,\eps}(e_k\otimes e_{-k}\otimes e_0)|^2\leq a |\lm_i|^{2k-2}\exp\left(-\frac{\pi}{2}|\lm_i|\right)\ ,
\end{eqnarray} for some explicit constant $a>0$. In fact, this is the actual order of the magnitude for the above value.

\subsubsection*{Remark} There is a natural trilinear functional on Whittaker models of representations of $G$. This is the model which appears in the Rankin-Selberg method as a result of unfolding. The above computation (and the similar one for Maass forms performed in \cite{BR-moscj}) shows that our normalization of the trilinear functional and the one coming from the Whittaker model coincide up to a constant of the absolute value one.

\subsubsection{Test vectors} Our construction is  very close to the construction we made in Section~\ref{consttest} for principal series representations, with appropriate modifications. We construct a test vector  $u_T(x,y)= v_T(x)\otimes v'_T(y)\in D_k\otimes D^*_k\subset \CH_{k-2}\otimes\CH_{-k}$ satisfying $H^{mod}_{k,\lm,\eps}(u_T)\geq \beta >0$ for $ T/2\leq |\lm|\leq T$ and  some constant $\beta>0$ independent of $\lm$. Vectors $v_T\in \CH_{k-2}$ and  $v'_T\in \CH_{-k}$ are first constructed  in the line model of induced representations, and then we relate these to vectors in the discrete series representation $D_k$.

 Choose a smooth non-negative function $\al\in C^\8(\br)$ with the support ${\rm supp}(\al)\subset [-0.1,0.1]$ and $\int_\br\al(x)dx=1$.  Consider the diagonal element $a_T={\rm diag}(T^{-\haf},T^\haf)\in G$. We define $w_T\in \CH_{k-2}$ by
\begin{eqnarray}\label{1-test-vect}
w_T=T^{\haf}\cdot\pi_{k-2}(a_T)\al\  .
                                        \end{eqnarray}
Recall that the action is given by $\pi_{k-2}({\rm diag}(a\inv,a)) v(x)=|a|^{2-k}v(a^2x)$.
We note that geometrically the vector $w_T$ is a small non-negative bump function around the point $0\in \br$, with the support in the interval $T\inv\cdot[-0.1,0.1]$, and satisfying $\int_{\br}w_T(x)dx=T^{\haf(1-k)}$. We now set $ v_T=q_k(w_T)\in D_k$. Note that $v_T=T^\haf\cdot \pi_{D_k}(a_T)\tilde v$ for some $\tilde v\in D_k$, i.e., the vector $v_T$ is obtained by the action of $G$ on some {\it fixed} vector in $D_k$. This will be crucial in what follows since we will need to estimate the supremum norm for the automorphic realization of the vector $v_T$.

The construction of the test vector  $v'_T\in D^*_k$ is slightly more complicated since we can not simply project a vector to $D^*_k\subset \CH_{-k}$ since the value of the functional $l^{ind}_{k,\lm,\eps}$ might change significantly. Let $\al$ be as above. We now view it as a vector in $\CH_{-k}$.  We choose  a smooth real valued function $\al'\in C^\8(\br)$ satisfying the following properties:
\begin{enumerate}
  \item ${\rm supp}(\al')\subset [M-0.1,M+0.1]$, where the parameter $M$ is to be chosen later,
  \item $\int_\br x^m\al'(x)dx=-\int_\br x^m\al(x)dx\ {\rm for}\ 0\leq m\leq k-2\ .$
\end{enumerate} The last condition implies that the vector $w=\al+\al'\in D_k^*$ since $\int_\br x^m w(x)dx=0$ for $0\leq m\leq k-2$.
We define now the second test vector by
\begin{eqnarray}\label{2-test-vect}
  v'_T=T^{\haf}\cdot\pi_{-k}(\left(
                                          \begin{smallmatrix}
1 & -1 \\
                                   &\ \ 1 \\
\end{smallmatrix}
                                        \right)a_T) w\ .
                                        \end{eqnarray}
Clearly we have $v'_T\in D^*_T$.

Recall that $\pi_{-k}({\rm diag}(a\inv,a)) v(x)=|a|^{k}v(a^2x)$, and $\pi_{-k}(\left(
                                          \begin{smallmatrix}
1 & -1 \\
                                   &\ \ 1 \\
\end{smallmatrix}
                                        \right)) v(x)=v(x+1)$.
Hence geometrically the vector $v'_T=\al_T+\al_T'$ is the sum of two bump functions $\al_T$ and $\al_T'$ with their supports satisfying
${\rm supp}(\al_T)\subset 1+T\inv[-0.1, 0.1]$ and ${\rm supp}(\al_T')\subset 1+T\inv[M-0.1, M+0.1]$, both near the point $1\in \br$. We also have $\al_T\geq 0$ and $\int_{\br}\al_T(x)dx=T^{\haf(k-1)}$.

We now set our test vector to $u_T=v_T\otimes v_T'\in D_k\otimes D^*_k$.  We want to show that for $ T/2\leq|\lm|\leq T$, \begin{eqnarray}\label{mod-lower-bnd-hol} |l^{mod}_{k,\lm,\eps}(u_T\otimes u)|\geq c'>0\end{eqnarray}  for some vector $u\in V_{\lm,\eps}$ with $||u||=1$, and with a constant $c'>0$ independent of $\lm$.

As we explained before, $l^{mod}_{k,\lm,\eps}(q_k(v)\otimes w\otimes u)= l^{ind}_{k,\lm,\eps}(v\otimes w\otimes u)$ for any triple $v\otimes w\otimes u\in \CH_{k-2}\otimes\CH_{-k}\otimes V_{\lm,\eps}$. Hence we work with $l^{ind}_{k,\lm,\eps}(w_T\otimes v'_T\otimes u)$ instead of $l^{mod}_{k,\lm,\eps}$. This value is given by an explicit integral. Let $K_{k,\lm}(x,y,z)$ be as in \eqref{K-kernel2}.
We have then
\begin{eqnarray}\label{triple-value-hol-test}
l^{ind}_{k,\lm,\eps}(w_T\otimes v'_T\otimes u)=\\ \int K_{k,\lm}(x,y,z)(sgn(x,y,z))^{\eps}\
 w_T(x)v'_T(y)u(z)\ dxdydz\ .\nonumber
\end{eqnarray} Hence it is enough to show that the absolute value of the integral
\begin{eqnarray}\label{triple-value-hol-z}
&&I_\lm(z)=\\ \int
 K_{k,\lm}(x,y,z) w_T(x)v'_T(y)\ dxdy&&=\langle K_{k,\lm}(x,y,z) , w_T(x)v'_T(y)\rangle \ \nonumber
\end{eqnarray} is not small for some interval of $z\in\br$. We have
$I_\lm(z)=K(z)+K'(z)$,  where we denote by $K(z)=\langle K_{k,\lm}(x,y,z) , w_T(x)\al_T(y)\rangle$ and by
$K'(z)=\langle K_{k,\lm}(x,y,z) , w_T(x)\al'_T(y)\rangle$.

Since the support of $w_T$ and of $\al_T$ is of the size $T\inv$, it is easy to see that $|K(z)|\geq c''>0$ for some constant $c''>0$, for $z\in [10,20]$ and for $|\lm|\leq T$. This is because the function $w_T(x)\al_T(y)$ is a non-negative function with the support in a small box of  size $T\inv$ around the point $(0,1)\in\br^2$, and that the gradient of
the function $K_{k,\lm}(x,y,z)$ is bounded by $|\lm|\leq T$ in this box.  Hence there are no significant cancellations in the integral   \eqref{triple-value-hol-z}. We normalized our vectors so that $\int_{\br^2}v_T(x)\al_T(y)dxdy=1$, and hence the integral $K(z)$ is not small for $z$ not near singularities of the kernel $K_{k,\lm}(x,y,z)$.   This part is identical to our argument in Section 5.3.4 of \cite{BR-moscj}.

We are left with the second term $K'(z)$. We want to show that there are no cancellations between two terms $K(z)$ and $K'(z)$ for  $z\in [10,20]$.  The reason for this is that while these terms are of about the same size, their arguments are different, and  not opposite if $\lm$ is not too small (e.g.,  $T/2\leq |\lm|\leq T$). Namely,  the argument of the kernel function $K_{k,\lm}$ in \eqref{triple-value-hol-z} on the support of $u_T$  is given by
\begin{eqnarray}\label{argument-K-lm}
\lm/2T\cdot \left(t_1(1-z\inv)-t_2(1+(z-1)\inv)\right)+ \\ \lm/2(\ln(z)+\ln(z-1)) +O(T\inv)\ ,\nonumber
\end{eqnarray} where $t_1\in [-0.1,0.1]$, $t_2\in [0.9, 1.1]$ for $\al_T$, and $t_2\in [M-0.1, M+0.1]$ for $\al_T'$.  By choosing appropriate value of $M$, we can see that the difference of these arguments  is not close to $0$ and $\pi$ for any {\it fixed}  $z\in [10,20]$ and $T/2\leq|\lm|\leq T$,  and hence integrals $K(z)$ and $K'(z)$ do not cancel each other since $u_T$ is real valued.

 Hence we have shown that $H_\lm(u_T)\geq c'>0$ for $T/2\leq |\lm|\leq T$, and some explicit $c'>0$ which is independent of $\lm$.

\subsubsection{Raising the level} We now discuss what happens when we change the level. Since our test vectors $v_T$ and $v_T'$ are {\it not} $K$-finite, we have to pass to automorphic functions on the space $X_\ell$. We use the standard notation  $j(g,z)=\det(g)^{-\haf} (cz+d)$ for
$g= \left(
                                          \begin{smallmatrix}
a & b \\
                               c    & d \\
\end{smallmatrix}
                                        \right)\in G^+$ and $z\in\uH$.
Let $\phi^k$ be a primitive holomorphic form of weight $k$ on $\uH$ for the subgroup $\G_0(N_0)$. We normalize $\phi^k$ by its norm on $Y_{N_0}$.  According to the well-known dictionary, we associate to $\phi^k$   the function $\phi_{e_k}\in C^\8(X_0)$ given by
\begin{eqnarray} \phi_{e_k}(g)=\phi^k(g(i))\cdot j(g,z)^{-k} \ , \end{eqnarray} where $z=g(i)$. In the opposite direction, we have
$\phi^k(g(i))=\phi_{e_k}(g)\cdot j(g,z)^{k}$.  We have the associated isometry $\nu_k=\nu_{\phi^k}:D_k\to C^\8(X_{N_0})$ which gives $\nu_k(e_k)=\phi_{e_k}$.

Let $\ell>1$ be an integer. We denote by $a_\ell={\rm diag}(
\ell^\haf ,\ell^{-\haf} )$.
For a given $\nu_k=\nu_{\phi^k}:D_k\to C^\8(X_{N_0})$, we construct the corresponding isometry ${\nu_k^\ell}:D_k\to C^\8(X_{\ell})$ as follows. For a vector $v\in D_k$, we consider the corresponding automorphic function
\begin{eqnarray}\label{nu-ell-nu-hol}
{^\ell\phi}_v(x)=\nu_k^\ell(v)(x)=\\ \nu_k(v)\left(a_\ell x\right)=\phi_v\left(a_\ell x\right)\  .\nonumber
\end{eqnarray} Obviously, we have $\sup_{X_\ell}|{^\ell\phi}_v|=\sup_{X_0}|{\phi}_v|$ for the same vector $v\in D_k$. We want to compare this to the classical normalization of old forms.
For the  lowest weight vector $e_k\in D_k$, we have with $z=g(i)$
 \begin{eqnarray}\label{nu-ell-nu-hol2}
{^\ell\phi}_{e_k}(g)=\phi_{e_k}
\left(a_\ell g\right)=\phi^k(\ell z)\cdot j(a_\ell g,i)^{-k}=\\ \ell^{\frac{k}{2}} \phi^k(\ell z)\cdot j(g,i)^{-k}\ . \nonumber
                                        \end{eqnarray}
On the other hand, classically, old forms are given by $ {^\ell\phi^k}( z)=\phi^k(\ell z)$. Hence we acquire the extra factor $ \ell^{\frac{k}{2}}$.

Since, as we noted, test  vectors $v_T\in D_k$ and $v'_T\in D^*_k$ are obtained by the (scaled) group action applied to  {\it fixed} vectors, and since the operation of raising the  level by $\ell$ does not change the supremum norm, we arrive at the following bound
\begin{eqnarray*}
\sup_{X_\ell}|{{^\ell\phi_{v_T}}(x)}|=\sup_{X_0}|{\phi_{v_T}(x)}|=\\ 
T^{\haf}\sup_{X_0}|{\phi_{\pi(a_T)w}(x)}|=
T^{\haf}\sup_{X_0}|{\phi_{w}(x)}|= \beta' T^\haf\ ,
\end{eqnarray*} for some constant $\beta'$. The same holds for the automorphic function ${^\ell\phi_{v'_T}}$. This implies that $||{^\ell\phi_{v_T}}{^\ell\phi_{v'_T}}||^2_{X_\ell}\leq \beta T^2$ for some $\beta>0$ independent of $\ell$ and $T$.

To summarize, we have proved the bound
\begin{eqnarray}\label{mean-value-hol}
\sum_{ \haf T\leq |\lm_i|\leq T}|\langle {^{\ell}}\phi^k\ \overline{^\ell\phi'^k},\phi_i\rangle_{Y_{\ell}}|^2\cdot\ell^{2k}|\lm_i|^{2-2k}e^{\haf\pi|\lm_i|}\leq A\cdot T^2\ ,
\end{eqnarray} where the summation is over an orthonormal basis of Maass forms for the subgroup $\G_0(\ell N_0)$.

The above argument also proves the case of forms with different level, i.e., the bound \eqref{mean-value-hol-dif-level}. Let  $\ell_1$ and $\ell_2$ be two co-prime integers.  Under our the normalization of measures on $X_\ell$, we see that for a vector $v\in V$, $L^2$-norms $||{^\ell\phi_v}||_{X_\ell}=||{^\ell\phi_v}||_{X_{\ell'}}$ are equal  if $\ell|\ell'$ (here we view the function ${^\ell\phi_v}$ as both $\G_0(\ell N_0)$-invariant function and as $\G_0(\ell' N_0)$-invariant function). Obviously, the supremum norms of $^\ell\phi_v$ on $X_\ell$ and on $X_{\ell'}$ are also coincide. Hence $$||{^{\ell_1}\phi_{v_T}}{^{\ell_2}\phi_{v'_T}}||^2_{X_{\ell_1\ell_2}}\leq  \haf||{^{\ell_1}\phi_{v_T}}^2||^2_{X_{\ell_1\ell_2}}+
\haf||{^{\ell_2}\phi_{v'_T}}||^2_{X_{\ell_1\ell_2}}\leq  \beta T^2\ .$$
This implies \eqref{mean-value-hol-dif-level}.
\qed


{\bf Acknowledgments.}  It is a pleasure to thank Joseph Bernstein for endless discussions concerning automorphic functions. I would like to thank Jeff Hoffstein for asking the question during the  Oberwolfach Workshop ``The Analytic Theory of Automorphic Forms" in August 2012 which led to this note, and the organizers of the workshop for their invitation.


\end{document}